

\input amstex
\documentstyle{amsppt}

\magnification=1200
\parskip 6pt
\pagewidth{5.4in}
\pageheight{7.1in}

\baselineskip=10pt
\expandafter\redefine\csname logo\string@\endcsname{}
\NoBlackBoxes                
\NoRunningHeads
\redefine\no{\noindent}

\define\C{\bold C}
\define\R{\bold R}
\define\Z{\bold Z}

\redefine\H{\bold H}

\define\al{\alpha}
\define\be{\beta} 
\define\ga{\gamma}
\define\de{\delta}
\define\la{\lambda}
\define\si{\sigma} 
\define\Si{\Sigma} 
\define\La{\Lambda}
\define\ep{\epsilon}
\define\Om{\Omega}

\define\om{\omega}

\define\Hol{\operatorname {Hol}}

\define\Harm{\operatorname {Harm}}
\define\alg{\operatorname{alg}}

\define\diag{\operatorname {diag}}

\redefine\dim{\operatorname {dim}}

\define\ad{\operatorname {ad}}

\define\sub{\subseteq}
\redefine\sup{\supseteq}

\redefine\b{\partial}
\define\pe{\perp}     

\def\pr{\prime} 
\define\st{\ \vert\ }   
\redefine\ll{\lq\lq}
\redefine\rr{\rq\rq\ }
\define\rrr{\rq\rq} 

\define\GC{G^{\C}}
\define\gc{\frak g^{\C}}

\define\grkcn{Gr_k(\C^n)}
\define\glnc{GL_n\C}
\define\frakglnc{\frak{gl}_n\C}

\redefine\F{F}
\redefine\H{H}

\define\LapGC{\La_+\GC}

\define\LaGC{\La\GC}
\define\LaalgGC{\La^{\alg}\GC}

\define\Grn{Gr^{(n)}}
\define\Hn{H^{(n)}}
\define\Ht{H^{(3)}}
\define\Grnalg{Gr^{(n)}_{\alg}}

\redefine\p{\frak p}
\redefine\n{\frak n}

\topmatter
\title An update on harmonic maps of finite uniton number, via
the zero curvature equation
\endtitle
\author Martin A. Guest
\endauthor
\endtopmatter

\document

The harmonic map equation is one amongst many examples of 
\ll integrable partial differential equations\rrr, that is,  p.d.e.
which can be written as \ll zero curvature equations\rrr.  It is a
relatively recent example, having been recognised as such (by
geometers) only in the 1980's.

Since there is, as yet, no universally accepted definition of integrability, the obvious
question that arises with any \ll integrable p.d.e.\rr is:  {\it how} integrable
is it? Classically, the ideal situation was a p.d.e. whose solutions can be represented
by concrete formulae involving elementary functions; of course this goal is
usually neither attainable (except in very special cases) nor desirable (concrete
formulae are often too complicated to be useful).  A more realistic
goal is the concept of a completely integrable Hamiltonian
system, but many interesting cases are not of this type either.

The purpose of this article is to review the equations for  {\it harmonic maps
of finite uniton number,} in a manner which demonstrates that these equations are
integrable in a very strong sense.  In fact the equations are integrable
in the most naive sense that all the solutions can be written down explicitly
in terms of \ll known\rr functions, but we shall not
make too much of this fact; we are more interested in
describing the solutions in ways likely to lead to new results and geometrical
insights.  It is this aspect which may be useful in other problems,
and which is indicative of the far-reaching
relations between integrable systems and geometry.

Since this particular equation is rather easy to describe (and solve), in
comparison with other \ll integrable\rr equations currently receiving the attention of
geometers, it provides a very instructive example.
Most of the previous work on harmonic maps of finite uniton number has
been published in conventional differential geometric contexts; we shall 
take a somewhat unorthodox point of view, which
may be more palatable to researchers in integrable systems theory. We
present this material in \S 1.

A secondary aim of this article
is to review the current state of affairs concerning
harmonic maps of finite uniton number.
The focus of recent articles in the differential geometry literature on
harmonic maps (from Riemann surfaces into Lie groups)
has been on maps which specifically do {\it not}
have finite uniton number; indeed,  harmonic
maps of finite uniton number are sometimes dismissed
on the vague grounds that such maps
have been \ll done\rrr, for example in \cite{Uh}. 
We would like to make the point that the latter paper, justly regarded as
a milestone in the theory of harmonic maps, was not in fact the end
of the story.  It would be more accurate to say that it was the beginning,
for its main achievement was to provoke geometers into
recognizing the zero curvature point of view, thus
establishing a very productive link with the world of integrable
systems.  For the convenience of the reader, in Appendix A we give a brief review of
the definitions and basic theory of harmonic maps of finite uniton number,
up to and including \cite{Uh}.

In \S 2, we give
several specific instances where this \ll integrable systems approach\rr
leads to concrete results which were either unattainable or left obscure by earlier
methods.  There are five topics, the first being the most fundamental:

\no (1) Canonical forms of complex extended solutions for
arbitrary compact Lie groups
\newline
(2) A Frenet frame construction for the case $G=U_n$
\newline
(3) Deformations of harmonic maps
\newline
(4) Topological properties of complex extended solutions
\newline
(5) Harmonic maps of finite uniton number into symmetric spaces

\no Of these, (1) and (2) show how to solve the equation for harmonic maps (of finite
uniton number) explicitly, in terms of \ll holomorphic data\rrr.
The method of (1) works for any compact Lie group $G$, but involves integrations
of meromorphic functions and hence troublesome residue considerations;
that of (2) involves only derivatives, but seems to work only for the case $G=U_n$.
Items (3) and (4) give applications to the structure of the space of
harmonic maps.  Since we restrict ourselves to maps into Lie groups
in (1)-(4), item (5) explains how these results extend to the case
of maps into symmetric spaces.  

To obtain these results, various refinements of the basic zero curvature
equation are necessary, as the pioneering concepts of \cite{Uh} 
(such as uniton factorizations) have turned out
to be awkward computational tools.  Some of our examples are previously
unpublished; in keeping with the expository style of the article
we postpone proofs of these new results to Appendices B and C.
We shall take the opportunity to indicate open problems along the way.

This article is a more detailed version of two lectures entitled 
\ll Harmonic maps: from
differential geometry to integrable systems\rr given at the National
Autonomous University of Mexico in September 2000, as part of the series
\ll Las Matem\acuteaccent aticas en el umbral del Siglo XXI\rrr. The author
thanks the other organisers of the 9-th MSJ-IRI for permitting him to submit
the article to these proceedings.  He is also grateful to Maarten Bergvelt,
Francis Burstall and Yoshihiro Ohnita for their collaboration on various
aspects of the work reviewed here, and to Josef Dorfmeister for several helpful
suggestions and for pointing out the references \cite{CFW}, \cite{Hn}.

\head
\S 1 Solving the harmonic map equation
\endhead

We shall consider harmonic maps $\phi:U\to G$, 
where $U$ is a simply connected open subset
of a compact Riemann surface $\Si$ and $G$ is a compact Lie group with a bi-invariant
Riemannian metric.  It will be convenient to start with the following (local)
formulation of the harmonic map equation:
$$
\H^{-1} \H^\pr=\frac 1\la \ \times \ \text{a meromorphic function $U\to\gc$}
\tag 1.1
$$
where the unknown function
$\H:U\times S^1\to\GC$ is assumed to be a meromorphic function
of  $z\in U$ and a smooth function of $\la\in S^1$ 
(the complex numbers of unit length).  We are writing $\GC$
for the complexification of $G$ here, and $\gc$ for its Lie algebra
(the complexification of the Lie algebra $\frak g$ of $G$).  
The notation $\H^\pr$ means
$\H_z$. Equation (1.1) imposes a strong condition on $\H$:
if we expand $\H$ as a Laurent series in $\la$ then $\H^{-1} \H^\pr$
will also have such an expansion, but (1.1) requires that all
coefficients except that of $\la^{-1}$ are zero.

In what sense is (1.1) a \ll zero curvature equation\rrr, and how is it
related to other definitions of harmonic maps?  
In general, if $\om$ is a connection form, its  curvature form is commonly written
as $d\om + \om\wedge\om$.  The zero curvature equation
is the equation $d\om + \om\wedge\om=0$. Locally, this condition
is equivalent
\footnote{
This equivalence is explained in detail in \cite{Sh}, where it is called
\ll the fundamental theorem of calculus\rrr.}
to the existence of a function $\H$ such that $\om=\H^{-1}d\H$.  
It is a general principle, supported by many examples, 
that if $\om$ is required to have some
particular form,  then the zero-curvature equation is equivalent to some
partial differential equation (for a function related to $\H$).
Equation (1.1) is of this type, and the p.d.e. in this case turns out to be the
well known harmonic map equation of differential geometry. We shall
explain this briefly in a moment (and in more detail in Appendix A).
Right from the start, however, we wish to emphasize that  
the above formulation has several advantages, the principal 
one being this:  {\it (1.1) is a system of first order
linear meromorphic ordinary differential equations.}  It is susceptible,
in principle, to classical methods of solution.  Moreover, it is a system depending
(in a very simple way) on the parameter $\la$, and in some cases this makes
solving the system very easy.

Consider, for definiteness, the equation $\H^{-1} \H^\pr=\frac 1\la A$ where
$A:U\to \GC$ is a given meromorphic function, with
$\GC=\glnc$  (all invertible $n\times n$ complex matrices).  Then
we seek a \ll fundamental solution matrix\rr $\H(z,\la)$ for
$$
(y_1^\pr(z),\dots,y_n^\pr(z))=\frac1\la (y_1(z),\dots,y_n(z)) A(z),
$$
i.e. a matrix whose rows are $n$ linearly independent solutions
of this system.

A formal solution of $\H^{-1} \H^\pr=\frac 1\la A$ of the form
$$
\H(z,\la)= \sum_{i\le 0} \H_i(z) {\la}^i
$$
may be obtained immediately by substitution, since the differential equation is
equivalent to the equations 
$$
\H_0^\pr = 0,\quad  \H_{-1}^\pr = \H_{0}A, \quad  \H_{-2}^\pr = \H_{-1}A, \dots
$$
which may be solved recursively for $\H_0,\H_{-1},\dots$.

Let us assume in addition that $A$ is nilpotent.  Then by taking $\H_{0}=I$ and all
constants of integration zero, we obtain a solution $\H(z,\la)$ which is
a polynomial in $\la^{-1}$.  In particular this formal solution 
converges, although it is not necessarily meromorphic because the integrations
may lead to logarithms.  {\it When such a solution is meromorphic, it corresponds
to a harmonic map $\Si\to G$ of finite uniton number, and we claim that all harmonic
maps $\Si\to G$ of finite uniton number are of this form.}  In other words, 
we claim that the
equation for harmonic maps $\Si\to G$ of finite uniton number is
equivalent to a system of differential equations which is solvable by quadrature ---
it is integrable in the most naive sense of the word.  This statement applies
to the case of any compact Lie group $G$, not just $G=U_n$, as we shall explain.

In principle the same method may be used to study arbitrary harmonic
maps (not just those of finite uniton number), 
but with formidable technical difficulties; the general situation
(which is far from understood) is reviewed in \cite{Do}.

\no{\it Example 1.2:}  Let us solve the equation $\H^{-1}\H^\pr=\frac1\la A$ for
$\H=\sum_{i\le0}\H_i\la^i$, where
$$
A=
\pmatrix
0 & u & v\\
0 & 0 & w\\
0 & 0 & 0 
\endpmatrix
$$
($u,v,w$ are given meromorphic functions of $z$).  Choosing $\H_0=I$ and all constants
of integration zero we obtain
$$
\H=I+\frac1\la 
\pmatrix
0 & \text{\eightpoint$\int u$} & \text{\eightpoint$\int v$}\\
0 & 0 & \text{\eightpoint$\int w$}\\
0 & 0 & 0
\endpmatrix
+\frac1{\la^2}
\pmatrix
{\ 0\ } & {\ 0\ }  & \text{\eightpoint$\int(w\int u)$}\\
{\ 0\ } & {\ 0\ } & 0\\
{\ 0\ } & {\ 0\ } & 0
\endpmatrix.\qed
$$

Keeping this example in mind, we return now to the general theory, to look at
the harmonic map equation in the context of the zero curvature equation
and some of the standard machinery of integrable systems.

We define a {\it complex extended solution} to be a map 
$\H:U\times S^1\to\GC$ which is holomorphic for $z\in U$ and smooth for
$\la\in S^1$, and is such that the Fourier series of $\H^{-1} \H^\pr$ 
is of the form
$$
\H^{-1} \H^\pr=\sum_{i\ge-1} A_i \la^i.
\tag 1.3
$$
In the language of Appendix A,  a complex extended solution is simply a
$\La\GC$-valued map which represents an extended solution $\F:U\to\Om G$, via
the identification $\Om G=\LaGC / \LapGC$.  

Let us briefly review this notation.  The based loop group $\Om G$ is the space
of (smooth) loops $\ga:S^1\to G$ such that $\ga(1)=e$, where $e$ is the identity
element of $G$. The complex loop group
$\LaGC$ is the space of all smooth maps $\ga:S^1\to \GC$, and $\LapGC$
is its subgroup consisting of maps such that $\ga$  extends
holomorphically to the unit disk. We refer
to \cite{Pr-Se} for further information on loop groups. However, the only fact
we really need in this section (and then only in order to explain
the relation with harmonic maps) is the above
identification and the equivalent statement that $\LaGC=\Om G\ \LapGC$ with $\Om
G\cap\LapGC=\{e\}$.  This implies
that any loop $\ga:S^1\to\GC$ may be factored uniquely as
$$
\ga=\ga_u\ga_+,\quad\quad \ga_u\in\Om G, \ \ga_+\in\LapGC.
$$
This factorization is a generalization of the Gram-Schmidt procedure
of linear algebra, the latter being equivalent to the factorization
of an invertible complex matrix $A$ in the form $A=A_uA_+$, where $A_u$ is
unitary and $A_+$ is upper triangular.

The connection between harmonic maps and complex
extended solutions is as follows.  If $\H$ is a complex extended solution, we regard
it as a map $\H:U\to\LaGC$, then define $\F=\H_u$, i.e. the
first factor of $\H=\H_u\H_+$.  The map $\phi(z)=\F(z,-1)$ is a harmonic
map from $U$ to $G$. Conversely, if $\phi:U\to G$ is harmonic, we obtain a corresponding
extended solution  $\F$ as in Appendix A, and hence
\footnote{To justify the passage from $\F$ to $\H$, some technical arguments are
needed, and these can be found in \cite{Do-Pe-Wu}.  In the case of maps
of finite uniton number, however, the existence of $\H$ is elementary.}
a complex extended solution $\H$.  

A complex extended solution $\H$ (or any corresponding
harmonic map) is said to have {\it finite uniton number} if $\H$
is a finite Laurent series in $\la$, i.e. a polynomial in $\la$ and
$\la^{-1}$. It is a result of \cite{Uh} (see Appendix A) that every
harmonic map $S^2\to G$ has finite uniton number.  By analogy with soliton
theory, harmonic maps of finite uniton number are sometimes called
unitons.

The concept of a complex extended solution originates from
several sources.  In \cite{Wa}, harmonic maps are identified with
certain holomorphic vector bundles, and complex extended solutions arise
as clutching functions.  (This vector bundle point
of view was developed further in \cite{An}.)  In \cite{Se}, harmonic maps are
identified with certain holomorphic maps into an infinite-dimensional Grassmannian,
and complex extended solutions arise (see \cite{Gu}) as representatives of such
maps with respect to natural coordinate charts.  In \cite{Do-Pe-Wu}, 
complex extended solutions arise in the above manner. They are often
called (holomorphic or meromorphic) extended frames.
All of these approaches aim to exploit the underlying complex geometry, a
point of view suggested by twistor theory.

Equation (1.3) (for complex extended solutions) has 
several theoretical advantages, in addition to the
undeniable practical advantage of being just a system of linear meromorphic ordinary
differential equations. We discuss
two such aspects next; both of them are fundamental
in the theory of integrable systems.

\no{\it (1) Wide availability of gauge transformations}

Although there is already some freedom in the choice of extended
solution for a given harmonic map, there is much greater freedom
in the choice of a complex extended solution:

\proclaim{Proposition 1.4} Let $\H$ be a complex extended solution.
Let $M$ be a $\GC$-valued map which is holomorphic in $z$ and
holomorphic in $\la$ for $0\le\vert\la\vert\le 1$.
Then the product $\H M$ is a complex extended solution.
\endproclaim

\demo{Proof}
It suffices to check that 
${(\H M)}^{-1} {(\H M)}^\pr = M^{-1} \H^{-1} \H^\pr M + M^{-1} M^\pr$ 
has (at worst) a simple pole at $\la=0$.  This is obvious, since
$\H^{-1} \H^\pr$ has at worst a simple pole at $\la=0$, and
$M$ is holomorphic there.
\qed\enddemo

The harmonic map associated to the complex extended solution $\H M$
is exactly the same as the one associated to $\H$, since $(\H M)_u = H_u$.
Therefore we gain the flexibility of choosing $M$ to suit our
purposes.  

\no{\it Example 1.5:} Let $\F$ be the extended solution associated to
$\H$, i.e. $\F=\H_u$ where $\H = \H_u \H_+$. Then we can say that
$\F$ is obtained by applying the gauge transformation
$\H\mapsto \H (\H_+)^{-1}=\F$.\qed

\no{\it Example 1.6:} Assume that there is a factorization
$\H = \H_- \H_+$ where $\H_-$, $\H_+$ (not
necessarily the same as in the previous example) are holomorphic 
for $0\le \vert\la\vert \le 1$, $1\le \vert\la\vert \le \infty$,
respectively.   
Then the gauge transformation $\H \mapsto \H (\H_+)^{-1}=\H_-$
produces a complex extended solution $\H_-$ such that
${(\H_-)}^{-1} {(\H_-)}^\pr$ is {\it linear} in $\la^{-1}$. (The latter
has a Laurent series with only non-positive powers of $\la$ by
definition, but the complex extended solution property means
that all terms except those in $\la^0$, $\la^{-1}$ are zero.)  
Now, it follows from the Birkhoff decomposition of
the loop group $\LaGC$ that such a factorization exists, perhaps
after translation of $\H$ by a constant loop, providing we
allow $\H_-$, $\H_+$ to have poles in $z$.  By multiplying
$\H_-$ on the right by $\H_-(z,\infty)^{-1}$, we can assume that
${(\H_-)}^{-1} {(\H_-)}^\pr$ is zero at $\la=\infty$, and hence
is of the form $\la^{-1}A(z)$, for some meromorphic map $A$. Thus we obtain a
complex extended solution having the special form of (1.1).   
This fundamental observation of \cite{Do-Pe-Wu} establishes the
equivalence of (1.1) and (1.3).\qed

\no{\it (2) Evident symmetry groups}

Symmetries are a prominent feature of integrable systems.  For complex
extended solutions, two symmetry groups arise very naturally:

\proclaim{Proposition 1.7} Let $\H$ be a complex extended solution of
finite uniton number.

\no (i) Let $\al\in\C^\ast$ (a nonzero complex number). Then the
map $\al\cdot\H(z,\la) = \H(z,\al\la)$ is a complex extended solution
of finite uniton number.

\no (ii) Let $\ga:S^1\to\GC$ be a map such that $\ga$ and $\ga^{-1}$ are finite
Laurent series in $\la$. Then the product $\ga\H$ is a complex extended solution
of finite uniton number.
\endproclaim

\demo{Proof}
(i) is obvious, and (ii) follows immediately from the formula
$(\ga\H)^{-1}(\ga\H)^\pr = \H^{-1}\H^\pr$.
\qed\enddemo

Thus we obtain actions of the group $\C^\ast$ and the \ll algebraic loop
group\rr $\LaalgGC$ (the set of maps $\ga$ which satisfy the
hypotheses of (ii) above) on the set of complex
extended solutions --- and hence on the set of harmonic maps ---
of finite uniton number.  While the effect on each of $\H$ and
$\H^{-1}\H^\pr$ is evident, the effect on the corresponding extended solution $\F=\H_u$
and the harmonic map $\phi(z)=\F(z,-1)$ is nontrivial --- unless $\ga\in \Om G$,
when the effect is simply to multiply $\F$ and $\phi$ by constants.
These actions coincide with the \ll circle action\rr and the
\ll dressing action\rr of \cite{Uh}, as was shown in \cite{Gu-Oh}.
The two actions do not commute, but they combine to give a natural
action of a semi-direct product group $\C^\ast\ltimes\LaalgGC$.

Having introduced the main properties of complex extended solutions,
we turn to the task of using them to study harmonic maps. The following
questions represent some reasonable goals:

\no{\it (Q1) For a given class of harmonic maps, how are such maps characterized
in terms of their complex extended solutions?}

Harmonic maps of finite uniton number were originally defined in \cite{Uh}
without reference to complex extended solutions; it follows from \cite{Uh} that
the original definition is equivalent to ours. Thus, our definition can be regarded
as a characterization of harmonic maps of finite uniton number in terms of their
complex extended solutions, i.e. that the complex extended
solutions are finite Laurent series in $\la$.  Another example is the class
of harmonic maps arising from the \ll twistor construction\rr  ---
the complex extended solutions of such maps can be characterized by a scaling condition
(see part (1) of \S 2).

\no{\it (Q2) Is there a canonical form amongst the complex extended solutions
associated to a given harmonic map?}

The large gauge freedom (and symmetry group) suggests the problem of
choosing a canonical representative, and hence parametrizing the space
of harmonic maps (or symmetry group orbits of harmonic maps) by specific meromorphic
functions. We shall in fact do this, later, for 
harmonic maps of finite uniton number, and it will be necessary
to make full use of both gauge transformations and
the symmetry group. For the moment, we give two
examples which provide motivation.   

\no{\it Example 1.8:} Let $f:U\to\C^n$ be a meromorphic function whose
derivatives $f,f^\pr,\dots,f^{(n-1)}$ are linearly independent at
almost all points of $U$. Let
$$
\H(z,\la) =
\pmatrix
\vert & \dots & \vert & \vert\\
f^{(n-1)} & \dots & f^\pr & f\\
\vert & \dots & \vert & \vert
\endpmatrix
\diag(\la^2,\dots,\la^2,\la,1,\dots,1)
$$
where $\diag(\la^2,\dots,\la^2,\la,1,\dots,1)$ is the $n\times n$
diagonal matrix with the first $n-i-1$ entries $\la^2$ and the last
$i$ entries $1$ (with $0\le i\le n-1$). Then $\H$ is a complex extended solution
associated to a harmonic map $\Si\to\C P^{n-1}$. It is well known
(see Appendix A) that all harmonic maps  $S^2\to\C P^{n-1}$ arise this way, so
the above formula for $\H$ may be regarded as a canonical form for such maps.\qed

\no{\it Example 1.9:} Let $p:U\to\C$ be a meromorphic function. Let
$$
\H(z,\la)= \exp \frac 1\la
\pmatrix
0 & p(z) \\
0 & 0
\endpmatrix
=
\pmatrix
1 & \frac 1\la p(z) \\
0 & 1
\endpmatrix
=
\pmatrix
1 & p(z) \\
0 & 1
\endpmatrix
\pmatrix
\frac 1\la & 0\\
0 & 1
\endpmatrix.
$$
Then $\H$ is a complex extended solution
associated to the holomorphic map $\Si\to\C P^1$
whose homogeneous coordinate expression is $z\mapsto [p(z);1]$.
This example is taken from \cite{Wa}.  It is equivalent to the case
$n=2$, $i=0$ of the previous example, and it represents a (different)
canonical form for such maps.\qed

\no{\it (Q3) How do geometrical properties of harmonic maps translate
into properties of complex extended solutions?}

For differential geometric properties in general (but primarily for
harmonic maps not of finite uniton number), we refer to \cite{Do} and
the references therein.  An instructive special case is that of
holomorphic maps $\Si\to S^2$, treated in \cite{Do-Pe-To}.
Other important special cases have been investigated in detail from the
point of view of differential geometry
in \cite{Br1}, \cite{Br2}, \cite{Br3}, \cite{Hn}, \cite{CFW}, \cite{Hs}.
In (3) and (4) of the next section we shall study a
topological property,  using complex extended solutions.

To demonstrate the effectiveness of complex extended solutions,
we conclude this section by
stating a very simple canonical form for complex extended solutions corresponding to
arbitrary harmonic maps $\Si\to U_n$ of finite uniton number 
(in particular for arbitrary harmonic maps $S^2\to U_n$).  In essence, this
is just a systematic generalization of Example 1.2.

Let $k\in\{0,1,\dots,n-1\}$ ($k$ will be the uniton number of the harmonic map,
as defined in Appendix A). Let
$v=(v_1,\dots,v_n)\in\Z^n$ with $k=v_1\ge v_2\ge \dots \ge v_n = 0$
and $v_i-v_{i+1}=0$ or $1$ for all $i$. We refer to $v$ as the \ll type\rr
of the harmonic map; it is in fact a Schubert symbol. Associated to $v$
there is a flag $\Cal F_v$ of subspaces of $\C^n$, and a parabolic
subgroup $P_v$ of $\glnc$, namely the group of all invertible
linear transformations of $\C^n$ which fix $\Cal F_v$.  The group $\glnc$
acts transitively on the space $\Om_v$ of all flags of type $v$, with isotropy
subgroup $P_v$ at $\Cal F_v$. It will be convenient to identify $\Cal F_v$
with the homomorphism $\ga_v:S^1\to U_n$, $\la\mapsto
\diag(\la^{v_1},\dots,\la^{v_n})$, and $\Om_v$ with the conjugacy class
of $\ga_v$. Let
the number of occurrences of $\la^k,\la^{k-1},\dots,1$ in
$\diag(\la^{v_1},\dots,\la^{v_n})$ be (respectively)
$a_0,a_1,\dots,a_k$.

Let $\p_v$ be the Lie algebra of $P_v$, and let $\p^0_v$ be the nil-radical
of $\p_v$; explicitly, this means that $\p_v$ consists of all complex matrices of the
form
$$
\pmatrix
A_{0,0} & A_{0,1} & \dots & A_{0,k-1} & A_{0,k} \\
0 & A_{1,1} & \dots & A_{1,k-1} & A_{1,k} \\
\vdots & \vdots & \ddots & \vdots & \vdots \\
0 & 0 & \dots & A_{k-1,k-1} & A_{k-1,k} \\
0 & 0 & \dots & 0 & A_{k,k}
\endpmatrix
$$
where $A_{i,j}$ is an $a_i\times a_j$ sub-matrix, and $\p^0_v$
consists of all complex matrices of the form
$$
\pmatrix
0 & A_{0,1} & \dots & A_{0,k-1} & A_{0,k} \\
0 &0 & \dots & A_{1,k-1} & A_{1,k} \\
\vdots & \vdots & \ddots & \vdots & \vdots \\
0 & 0 & \dots & 0 & A_{k-1,k} \\
0 & 0 & \dots & 0 & 0
\endpmatrix.
$$
We have the descending central series of $\p^0_v$,
$$
\p^0_v \sup \p^1_v \sup \dots \sup \p^{k-1}_v \sup \p^k_v = \{0\},
$$
defined by $\p^{i}_v=[\p^{i-1}_v,\p^{0}_v]$. The Lie algebra $\p^{i}_v$
consists of all $n\times n$ complex matrices in \ll block form\rr
$(A_{\al,\be})_{1\le\al,\be\le k}$, with $A_{\al,\be}=0$ for
$\al\ge \be-i$.

For each $i\in\{1,\dots,k\}$, let $B_i$ be a $\p^{i-1}_v$-valued meromorphic function.
Consider the map
$$
\H(z,\la)=\exp\,B(z,\la)\quad\text{where}\quad
B(z,\la)=\frac1\la B_1(z)+\frac1{\la^2} B_2(z) +\dots + \frac1{\la^{k}}B_{k}(z).
$$
Then ${\H}^{-1}\H^\pr$ is a polynomial in $\la^{-1}$ with no constant
term, and $\H$ is a complex extended solution if
and only if the only nonzero coefficient in this 
polynomial is that of $\la^{-1}$ itself. 
The coefficient of $\la^{-1}$ is evidently $B_1^\pr$, so we can say that
$\H$ is a complex extended solution if and only if
${\H}^{-1}\H^\pr=\frac1\la B_1^\pr(z)$.  By the well known
formula for the derivative of the exponential map (\cite{He}, Chapter 2, Theorem 1.7),
the complex extended solution condition is therefore
$$
B^\pr-\frac{1}{2!}(\ad B)B^\pr+\frac{1}{3!}(\ad B)^2B^\pr-
\frac{1}{4!}(\ad B)^3B^\pr+\dots \quad = \quad \frac1\la B_1^\pr
$$
where $\ad B X$ means $BX-XB$.
This condition is --- as expected --- a system of meromorphic ordinary
differential equations for $B_1,\dots,B_{k}$ which can be integrated recursively.
Indeed, for $i=2,\dots,k$, the equation arising from the coefficient of $\la^i$ expresses
$B_i^\pr$ as a polynomial in terms of $B_1,\dots,B_{i-1}$ and
$B_1^\pr,\dots,B_{i-1}^\pr$.  For example, the first equation, for $i=2$, is
$B_2^\pr = \frac12(B_1B_1^\pr - B_1^\pr B_1)$.

\no{\it Example 1.10:} For the group $U_3$ and the homomorphism
$\ga_v(\la)=\diag(\la^2,\la,1)$, $B$ is of the form
$$
B=\frac 1\la B_1+\frac1{\la^2} B_2
=
\frac 1\la
\pmatrix
0 & a & b\\
0 & 0 & c\\
0 & 0 & 0
\endpmatrix
+
\frac 1{\la^2}
\pmatrix
0 & 0 & d\\
0 & 0 & 0 \\
0 & 0 & 0 
\endpmatrix
$$
where $a,b,c,d$ are meromorphic functions of $z$.  
The complex extended solution condition
is $B_2^\pr = \frac12(B_1B_1^\pr - B_1^\pr B_1)$, which is just
$$
d^\pr = \frac12(ac^\pr-a^\pr c).
$$
This can be solved directly, by choosing $a,b,c$ arbitrarily and then integrating
to obtain $d$. To compare this complex extended solution with the one in Example 1.2,
we compute
$$
\H=\exp\,B
=
\pmatrix
1&0&0\\
0&1&0\\
0&0&1
\endpmatrix
+
\frac 1\la
\pmatrix
0 & a & b\\
0 & 0 & c\\
0 & 0 & 0
\endpmatrix
+
\frac 1{\la^2}
\pmatrix
0 & 0 & d+\frac12 ac\\
0 & 0 & 0 \\
0 & 0 & 0 
\endpmatrix.
$$
This is exactly the formula of Example 1.2, with
$u=a^\pr$, $v=b^\pr$, $w=c^\pr$.  (Note that the derivative
of $d+\frac12ac$ is then $d^\pr +\frac 12(a^\pr c + a c^\pr) = 
\frac12(ac^\pr-a^\pr c)+\frac 12(a^\pr c + a c^\pr)
= ac^\pr = w\int u$.)\qed

It turns out that all harmonic maps $\Si\to U_n$ 
of finite uniton number arise this way:

\proclaim{Theorem 1.11} (i) Let $\H$ be a complex extended solution, corresponding to
a harmonic map $\phi:\Si\to U_n$ of finite uniton number. Then
there exists a gauge transformation (in the sense of Proposition 1.4)
which converts $\H$ to the above canonical form $\exp\,B$, for some $v$ and
some meromorphic $B_1,\dots,B_k$.

\no(ii) Conversely, let $B_1:\Si\to\p^0_v$ be any meromorphic function.
Let $B_2,\dots,B_k$ be functions obtained by solving recursively
the system of equations obtained from equating coefficients of powers of $\la$
in
$$
B^\pr-\frac{1}{2!}(\ad B)B^\pr+\frac{1}{3!}(\ad B)^2B^\pr-
\frac{1}{4!}(\ad B)^3B^\pr+\dots \quad = \quad \frac1\la B_1^\pr.
$$
If $B_2,\dots,B_k$ are meromorphic, then $B$
defines a complex extended solution $\H=\exp\,B$ (corresponding to
a harmonic map $\phi:\Si\to U_n$ of finite uniton number).

\no(iii) The effect of changing the constants of integration in (ii)
is to change $\H$ by the (dressing) action of the group $\LaalgGC$
in the sense of Proposition 1.7 (ii)).
\endproclaim

\demo{Proof} (i) is a special case of the results of \cite{Bu-Gu}; we
shall review this in the next section.  (ii) is immediate from the
construction above. (iii) follows from two facts. First, the dressing
action of Proposition 1.7 (ii) does not change $\H^{-1}\H^\pr$,
hence only the constants of integration are affected. Second, $\H^{-1}\H^\pr$
determines $\H$ up to multiplication on the left by a map $\ga:S^1\to \GC$,
i.e. up to the dressing action.
\qed\enddemo

The above construction may be regarded as an Ansatz which, when
used in (1.1),  happens to produce all harmonic maps of 
finite uniton number.  We shall explain (geometrically) 
why it works in the next section.
From the point of view of integrable systems, it is notable
that the canonical form (with all constants of integration taken
to be zero) can be achieved by a combination of gauge transformations and
dressing transformations.

It is important to note that the passage from the
complex extended solution $\H$ to 
the corresponding extended solution $\F=\H_u$ and harmonic map 
$\phi=\F\vert_{\la=-1}$ is a purely (real) algebraic operation.  
This is because the factorization $\ga=\ga_u\ga_+$ of an
algebraic loop $\ga$ is a purely algebraic operation.  Geometrically,
the factorization is an infinite-dimensional version of the Gram-Schmidt
orthogonalization procedure, but when $\ga$ is algebraic the
procedure takes place in a finite-dimensional vector space (this can be seen
from the Grassmannian model of the loop group $\Om U_n$).

Finally, we remark that the construction shows that the number
of meromorphic functions $\Si\to\C$ required to produce a harmonic
map $\Si\to U_n$ of type $v$ is (at most) the
complex dimension of the vector space $\p^0_v$.  We have therefore
obtained a very precise and complete solution of the harmonic
map equation in this situation.  The only blemish
is the fact that {\it arbitrary} meromorphic
functions are not allowed, only those which lead to meromorphic
functions in the integration procedure described above.  This problem
will be addressed in (2) of the next section.

\head
\S 2 Explanations, applications, generalizations
\endhead

\no{\it (1) Canonical forms of complex extended solutions for
arbitrary compact Lie groups}

Theorem 1.11 has a rather surprising origin; surprising, that is,
from the point of view of differential geometry, but one that is
typical of the integrable systems approach.
This is explained in
\cite{Bu-Gu}, where a canonical form was given for complex extended solutions
$\H:\Si\to\GC$ of finite uniton number, for any compact Lie group $G$.  
To avoid introducing further Lie algebraic notation, 
we shall just discuss the case $G=U_n$ here.

The key ingredient is the loop group $\Om U_n$. This has a well known
Morse-theoretic decomposition into stable manifolds with respect to the
energy functional $E:\Om U_n\to\R$, $\ga\mapsto \int_{S^1} \vert d\ga \vert^2$.
The critical points of $E$ are the geodesic loops, i.e. the homomorpisms.
The connected critical manifolds are the conjugacy classes $\Om_v$ of the homomorphisms
$\ga_v$ with $v=\diag(\la^{v_1},\dots,\la^{v_n})$ and $v_1\ge\dots\ge v_n$.
Although $\Om U_n$ and each stable manifold is infinite-dimensional, it 
turns out that each unstable manifold $U_v$ is finite-dimensional
(with the structure of a vector bundle of rank $r(v)$ over $\Om_v$,
where $r(v)$ is the index of the geodesic $\ga_v$). The union of these unstable
manifolds is known (\cite{Pr}) to be the algebraic loop group
$\Om^{\alg}U_n$, a proper subset of $\Om U_n$.  Without loss
of generality (\cite{Uh}, \cite{Se}), for the purpose of studying harmonic maps
of finite uniton number, it suffices to consider extended solutions
$\F:\Si\to \Om^{\alg}U_n$.  

Since $\F$ is holomorphic,
its image must be contained in the closure of a single unstable manifold
$U_v$, for some particular $v$.  The main result of \cite{Bu-Gu} is that
it suffices to consider only a {\it finite number} of such $v$, namely those
satisfying the additional conditions $v_i-v_{i+1}=0$ or $1$ for all $i$,
and $v_n=0$. A proof of this fact was given in Chapter 20 of
\cite{Gu}, using the Grassmannian model of $\Om U_n$, whereas the
proof of \cite{Bu-Gu} uses the loop group $\Om G$ directly, and works for
any $G$.

Since each $v_i-v_{i+1}$ can be $0$ or $1$, 
there are $2^{n-1}$ possible \ll types\rr of harmonic maps, 
and in fact they all occur, i.e. no further reduction is possible.  
For $n=3$ the types are $(2,1,0)$, $(1,1,0)$, $(1,0,0)$, and $(0,0,0)$. 
The first case gives the harmonic maps of Example 1.10. The second and third types
correspond to holomorphic maps into $Gr_2(\C^3)$ and $\C P^2$, respectively.  The fourth
type corresponds to constant maps.

Now, the structure of the unstable manifold $U_v$ is well known, and in
particular it has a \ll big cell\rr consisting of loops of the form
$$
[\exp\,(P_0+\la P_1 +\dots + \la^{k-1}P_{k-1})]\ga_v
$$
where $P_i\in \p^i_v$.  
The canonical form arises by expressing the complex
extended solution equation in the \ll new coordinates\rr 
$P_0,\dots,P_{k-1}$, as follows.
For each $i\in\{0,1,\dots,k-1\}$, let $C_i$ be a $\p^i_v$-valued meromorphic
function. Consider the map
$$
\H(z,\la)=[\exp\,C(z,\la)]\ga_v(\la)\quad\text{where}\quad
C(z,\la)=C_0(z)+\la C_1(z) +\dots + \la^{k-1}C_{k-1}(z).
$$
Using the formula for the derivative of the exponential map,
it is easy to verify that $\H$
is a complex extended solution if and only if, 
for each $i=0,1,\dots,k-2$, the coefficient of $\la^i$ in 
$$
\sum_{n\ge 0} \frac{(-1)^n}{(n+1)!} (\ad\,C)^n C^\pr
$$
has zero component in $\p^{i+1}_v$.  We obtain a system of 
meromorphic differential equations which can be solved recursively for
$C_0,\dots,C_{k-1}$. So far this appears somewhat different from the canonical
form $\exp\, B$ given at the end of the last section; however, one has
$$
\exp\,B = {\ga_v}^{-1} \exp\, C \ga_v = \exp\, {\ga_v}^{-1} C\ga_v,
$$
i.e. the two versions differ only by a \ll trivial\rr dressing
transformation.  

All this is valid for maps into a general compact Lie group $G$.
More precisely, we can say that Theorem 1.11 is valid for such groups, providing
the Lie algebras $\p^i_v$ are defined Lie algebraically as in \cite{Bu-Gu}.
An immediate consequence is an upper bound on
the (minimal) uniton number.  For $G=U_n$
it is $n-1$; for general $G$ the upper bound can be expressed root-theoretically, and
the results for the simple groups are as follows:

\bigpagebreak

\settabs 4\columns

\+  & $G$ & $\text{upper bound on uniton number}$ &\cr

\+ & ------ & ----------------------------------------- &\cr

\+ & $SU_n$ & $n-1$ &\cr

\+ & $SO_{2n+1}$ & $2n-1$ &\cr

\+ & $Sp_n$ & $2n-1$ &\cr

\+ & $SO_{2n}$ & $2n-3$ &\cr

\+ & $G_2$ & $5$ &\cr

\+ & $F_4$ & $11$ &\cr

\+ & $E_6$ & $11$ &\cr

\+ & $E_7$ & $17$ &\cr

\+ & $E_8$ & $29$ &\cr

\medpagebreak

Another consequence of the canonical form (and the algebraic nature
\footnote{The algebraic nature of the factorization can be seen from the
Grassmannian model in the case $G=U_n$; in the general case it follows by
taking a faithful unitary representation of the compact Lie group $G$.}
of the Iwasawa factorization) is the fact that all harmonic maps $\Si\to G$ of
finite uniton number are (real) algebraic functions of meromorphic
functions on $\Si$. This answers a question raised at the end of \cite{Wo}, where 
a construction of harmonic maps $S^2\to U_n$ was distilled from
the uniton factorization of \cite{Uh}.

It should be noted, incidentally, 
that the above construction provides a very satisfactory
implementation of the \ll DPW method\rr (introduced in \cite{Do-Pe-Wu} and
surveyed in \cite{Do}) for the case of harmonic maps of finite uniton number.
The DPW meromorphic potential (now usually called the normalized potential) 
is simply $B_1^\pr dz$.

There are special families
of harmonic maps of finite uniton number, including those sometimes
known as superminimal or isotropic, which have been given a general formulation
via the \ll twistor construction\rr in \cite{Bu-Ra}, \cite{Br3} and \cite{Sa}. 
These maps have
a simple characterization in terms of extended solutions:  they correspond to
extended solutions $\F$ which are $S^1$-invariant, where $S^1$ acts on the
loop group $\Om G$ by $\al\cdot\ga(\la)=\ga(\al\la)\ga(\al)^{-1}$. Geometrically,
this means that the image of $\F$ lies in a conjugacy class of homomorphisms 
(a connected component of the set of critical points of the function  
$\ga\mapsto\int_{S^1}\vert d\ga\vert^2$), i.e. a generalized flag manifold
of the group $G$.
In the case $G=U_n$, it means that $\F$ takes values in $\Om_v$ for some $v$.
The corresponding harmonic map $\phi$ always takes values in a symmetric space
of $G$ (embedded totally geodesically in $G$).  The relation between
$\F$ and $\phi$ is particularly simple in this situation:  $\phi$ is
obtained by composing $\F$ with a natural projection from the generalized
flag manifold to the symmetric space. Such projection maps are
called (generalized) twistor fibrations. (The twistor fibration $\C P^3\to S^4$
of Yang-Mills theory is an example; here $G=SO_5$ or $Sp_2$.)  Harmonic
maps of this type were originally studied without reference to loop groups,
although the loop group formulation provides a unifying framework.

Let us consider now the harmonic maps whose extended solutions
take values in
the conjugacy class $\Om$ of a particular homomorphism $\ga\in\Om G$.
The complex extended solutions $\H$ of such maps
are characterized by the following scaling condition:
$$
\H(z,\al\la)=\H(z,\la)\ga(\al)\quad\text{for all $\al\in S^1$}.
$$
This fact is an immediate consequence of the canonical form
$\H(z,\la)=[\exp\,C(z)]\ga(\la)$.
Here $C$ is a meromorphic function taking values in the nil-radical
of the parabolic subalgebra determined by $\ga$, i.e. in the
\ll big cell\rr of the generalized flag manifold $\Om$.
(In other words, we have  $C_1=\dots=C_{k-1}=0$ and $C=C_0$.)
A feature of this scaling condition
is that it distinguishes the particular twistor fibration.
The dressing action of Proposition 1.7 (ii) preserves the
generalized flag manifold $\Om$; in fact it reduces to the
natural action of $\GC$ on this space.

The complex extended solution equation
for harmonic maps of this special type can be integrated directly, exactly
as in the case of general maps of finite uniton number. Indeed, this
was done by Bryant in \cite{Br3}, 
long before the introduction of extended solutions,
at least in the special case of uniton number $2$.  Several authors have studied
the meromorphic o.d.e. system for harmonic maps into spheres, including
Bryant himself in \cite{Br1}, \cite{Br2} and then \cite{Hn}, \cite{CFW}, \cite{Hs}.

\no{\it (2) A Frenet frame construction} 
\footnote{The results in (2) have been
obtained in joint work with Francis Burstall.}
{\it for the case $G=U_n$}

The canonical form of (1), together with the \ll integrable\rr
meromorphic o.d.e. for its coefficient functions, shows that
harmonic maps $\Si\to G$ of finite uniton number are parametrized by 
collections of meromorphic functions on $\Si$.  However, it is not
clear from this description which meromorphic functions are
disallowed on the grounds that they give rise to logarithms during
the integration process.

It turns out that, in the case $G=U_n$, this difficulty may be avoided: 
there is an alternative parametrization of (complex extended solutions
corresponding to) harmonic maps $\Si\to U_n$ of finite uniton number, 
in which the initial
data consists of arbitrary meromorphic functions and their derivatives; no
integrations are needed.  One way to obtain this new parametrization is 
to make a change of variable in the
meromorphic o.d.e. for complex extended solutions:

\no{\it Example 2.1:} Consider the o.d.e. 
$d^\pr - \frac12(ac^\pr-a^\pr c)=0$ for complex extended solutions
of the type $G=U_3$, $v=(2,1,0)$ (Example 1.10).  Instead of choosing
$a,b,c$ and integrating to obtain $d$, let us introduce new variables
$$
\al=d+\frac12 ac,\quad
\be = c,\quad
\ga=a,\quad
\de=b.
$$
Then the o.d.e. simplifies to $\al^\pr=\be^\pr \ga$.  {\it This may be
solved without integration:}  choose meromorphic functions
$\al,\be,\de$, then obtain $\ga$ as $\ga=\al^\pr/\be^\pr$. Thus, all
harmonic maps of this type may be constructed by choosing as initial data
three meromorphic functions $\al,\be,\de$ and then performing a series
of derivatives and algebraic operations.\qed

It can be shown that the same phenomenon occurs for complex
extended solutions of type $v$ for $U_n$, in general.
We shall not prove this here, because there is a more geometrical, 
though less economical, version of the parametrization,  
in which the initial data consists of a collection of meromorphic maps
$\Si\to\C^{n}$.  This is in the spirit of the Eells-Wood description of harmonic
maps into $\C P^{n-1}$ (see Theorem A.1 of Appendix A), 
and our result can be viewed as
the natural generalization to the case where the target manifold is $U_n$.
A similar result holds for $\grkcn$.

The appropriate context for this version --- and this is why our argument
is restricted to the case $G=U_n$ --- is the Grassmannian model
$\Grn=\Om U_n$. As explained in Appendix A, the complex extended solution
equation for a holomorphic map $W:\Si\to\Grn$ is
$$
\la W^\pr\sub W. \tag 2.2
$$
For each $z\in\Si$, $W(z)$ is a linear subspace of the Hilbert space
$$
\Hn=L^2(S^1,\C^n)=\bigoplus_{i\in\Z} \la^i \C^n.
$$
However, for maps of uniton number $k$, it suffices (by \cite{Se}) to consider the
case
$$
\la^k\Hn_+ \sub W \sub \Hn_+
$$
where
$$
\Hn_+=\bigoplus_{i\ge 0} \la^i \C^n.
$$
Thus, $W$ is a holomorphic map to a {\it finite-dimensional} Grassmannian manifold.
Although this means that the original Hilbert space 
$\Hn$ will play no essential role, we
continue to use it for notational convenience; the reader should bear in mind
that we always work in the finite-dimensional vector space
$\Hn_+/\la^k\Hn_+ \cong \oplus_{i=0}^{k-1} \la^i \C^n \cong \C^{kn}$.
If $W$ is spanned locally by $\C^{kn}$-valued holomorphic functions 
$s_1,\dots,s_r$ (mod $\la^k\Hn_+$), we write $W=[s_1]+\dots+[s_r]+\la^k\Hn_+$. 
With this convention, the notation $W^\pr$ means 
$$
W^\pr = [s_1]^\pr +\dots+[s_r]^\pr 
 = [s_1]+\dots+[s_r]+[s^\pr_1]+\dots+[s^\pr_r] +\la^k\Hn_+.
$$

Equation (2.2) imposes a nontrivial condition on the holomorphic
map $W$, but  all such maps may be generated
from arbitrary holomorphic maps in a simple way:

\proclaim{Proposition 2.3} Let $X:\Si\to Gr_s(\oplus_{i=0}^{k-1} \la^i \C^n)$
be any holomorphic map.  Define $W$ by
$$
W=X+\la X^\pr + \la^2 X^{\pr\pr} + \dots + \la^{k-1} X^{(k-1)} + \la^k\Hn_+.
$$
Then $W$ is a solution of (2.2), and all solutions 
$W:\Si\to \Hn_+/\la^k\Hn_+$ of (2.2) arise this way.
\endproclaim

\demo{Proof} Given $X$, it is obvious that $W$ satisfies (2.2). Conversely, given
any solution $W$ of (2.2), we may take any $X$ such that $W=X+\la^k\Hn_+$.
\qed\enddemo

As it stands, this result is quite trivial, although it demonstrates already that
harmonic maps of finite uniton number may be constructed from \ll unconstrained\rr
holomorphic data using only derivatives and algebraic operations.  It is of
more interest to know how to choose  holomorphic data
for a {\it particular type} of solution, perhaps. 
As a concrete example, we shall do this below 
for the cases $G=U_3$ and $U_4$.

We need the following additional notation.
Any holomorphic map $X:\Si\to Gr_s(\oplus_{i=0}^{k-1} \la^i \C^n)$
may be written in the form
$$
X=X_0+\la X_1 + \la^2 X_2 + \dots + \la^{k-1} X_{k-1},
$$
where $X_i:\Si\to Gr_{s_i}(\oplus_{j=0}^{k-i-1} \la^j \C^n)$
for some $s_i$; we may assume that the image of $X_i$ is not contained
in $Gr_{s_i}(\oplus_{j=0}^{k-i-2} \la^j \C^n)$.  (In other words, we
take the echelon form of $X$ with respect to the flag given by
the subspaces $\oplus_{j=0}^{k-i-1} \la^j \C^n$.)

We shall specify suitable maps 
$X_0,X_1,\dots$ for each type of (non-constant) harmonic map
of finite uniton number into $U_3$ and $U_4$ in the tables below. In each
case, $l,m,\dots:\Si\to\C^n$ are arbitrary vector-valued meromorphic functions.
We regard these functions as \ll Frenet frame data\rrr, analogous to the
$\C^n$-valued meromorphic function $f$ in Theorem A.1 of Appendix A.

\bigpagebreak

\settabs 4\columns
\+   $\text{Frenet frame data for maps $\Si\to U_3$:}$\cr
\+   \cr
\+  $(2,1,0)$ & $X_0 = [l+\la m]$ & $X_1=0$ \cr
\+  $(1,1,0)$ & $X_0 = [l]$ & \cr
\+   $(1,0,0)$ & $X_0 = [l]+[m]$ & \cr

\bigpagebreak

\settabs 4\columns
\+   $\text{Frenet frame data for maps $\Si\to U_4$:}$ \cr
\+   \cr
\+  $(3,2,1,0)$ & $X_0 = [l+\la m +\la^2 n]$ & \quad\quad\quad$X_1=0$ & \quad$X_2=0$\cr
\+   $(2,2,1,0)$ & $X_0 = [l+\la m]$ & \quad\quad\quad$ X_1=0$  \cr
\+   $(2,1,1,0)$ &$X_0 = [l+\la m]$ & \quad\quad\quad$X_1=[n]$  \cr
\+ $(2,1,0,0)$ & $X_0 = [l+\la m]+[l^\pr+\la n]$ & \quad\quad\quad$X_1=0$ \cr
\+ $(1,1,1,0)$ & $X_0 = [l]$ & \cr
\+ $(1,1,0,0)$ & $X_0 = [l]+[m]$ & \cr
\+ $(1,0,0,0)$ & $X_0 = [l]+[m]+[n]$ & \cr

\medpagebreak

\comment

\no$(2,1,0)$\quad $X_0 = [l+\la m], X_1=0$

\no$(1,1,0)$\quad $X_0 = [l]$

\no$(1,0,0)$\quad $X_0 = [l]+[m]$

\medpagebreak

\no$(3,2,1,0)$\quad $X_0 = [l+\la m +\la^2 n], X_1=0, X_2=0$

\no$(2,2,1,0)$\quad $X_0 = [l+\la m], X_1=0$

\no$(2,1,1,0)$\quad $X_0 = [l+\la m], X_1=[n]$

\no$(2,1,0,0)$\quad $X_0 = [l+\la m]+[l^\pr+\la n], X_1=0$

\no$(1,1,1,0)$\quad $X_0 = [l]$

\no$(1,1,0,0)$\quad $X_0 = [l]+[m]$

\no$(1,0,0,0)$\quad $X_0 = [l]+[m]+[n]$

\endcomment

\proclaim{Theorem 2.4}  All harmonic maps $\Si\to U_n$
of finite uniton number, for $n=3$ or $4$,
arise through the above construction.  That is, for any choice of 
$\C^n$-valued meromorphic functions $l,m,\dots$, the map 
$=X+\la X^\pr + \la^2 X^{\pr\pr} + \dots + \la^{k-1} X^{(k-1)} + \la^k\Hn_+$
(where $X=X_0+\la X_1 + \la^2 X_2 + \dots + \la^{k-1} X_{k-1}$)
is an extended solution corresponding to a harmonic map $\Si\to U_n$
of finite uniton number, and all such maps arise this way.
\endproclaim

We relegate the proof, which is an explicit computation
based on the canonical form given in part (1), to Appendix B.  
The proof shows that {\it generic} choices of $l,m,\dots$ 
produce extended solutions of the indicated type, while for
certain special choices an extended solution of simpler type may be
obtained.  It also shows that $l,m,\dots$ are expressable in terms
of the data of the canonical form via a change of variable like that in
Example 2.1.  

The same method
works for $U_n$ in general.  As the definitions of $X_0,X_1,\dots$
in the general case are increasingly complicated (and non-canonical), 
we just remark that the number of $\C^n$-valued
meromorphic functions $l,m,\dots$ needed to construct a harmonic
map $\Si\to U_n$ of finite uniton number is at most $n-1$.  To see this, let
the number of occurrences of $\la^k,\la^{k-1},\dots,1$ in
$\ga_v(\la)=\diag(\la^{v_1},\dots,\la^{v_n})$ be (respectively)
$a_0,a_1,\dots,a_k$.  An upper bound on the number of
$\C^n$-valued functions required can be computed by considering the
\ll worst case\rr $a_k\ge a_{k-1}\ge\dots \ge a_0$.  The number in this case is
at most
$$
a_0 k + (a_1-a_0)(k-1) + \dots + (a_k-a_{k-1})0 = 
a_0+a_1+\dots+a_{k-1} = n-a_k \le n-1
$$
as claimed.  Whether the method can be extended to other compact Lie
groups $G$ is an open problem.

Let us try to make an honest assessment of the
various methods introduced so far.  First, the  method of (1)
gives a solution in terms of the minimum number of meromorphic
functions $\Si\to\C$, i.e. the nonzero component functions of $B_1$,
and this holds for any compact Lie group $G$.
However,  one has to exclude meromorphic
functions which produce nonzero residues in the integration process
(and it does not seem to be easy to characterize such functions).
For $G=U_n$, the change of variables avoids this problem, and allows
us to express the solution in terms of unconstrained holomorphic data,
but specifying the change of variable needed
for each type of harmonic map is complicated.  Finally, the method of
Theorem 2.4 gives a more systematic solution using unconstrained initial data,
again in the case $G=U_n$,
but this data is certainly not minimal. For example, eight meromorphic
functions (the components of $l$ and $m$)
are needed for maps into $U_3$ of type $(2,1,0)$,
whereas we know from Example 2.1 that four are enough.  This kind of problem
was visible already in Examples 1.8 and 1.9, in the description of
holomorphic maps $S^2\to S^2$. 
It seems, in conclusion, that what is the \ll best\rr method depends on
what one is trying to do.

\no{\it (3) Deformations of harmonic maps}

The dressing action is an action of an infinite-dimensional
Lie group on the space of harmonic maps
(or solutions of any integrable p.d.e.), and it has important theoretical
implications.  From a purely practical point of view it produces many new
solutions from a given solution, and it suggests the problem of describing
the orbits of the action.  In the case of harmonic maps of finite
uniton number (in contrast to examples like the KdV equation),
neither of these leads to significant new results.
The reason is that the orbits of the dressing action are too small (and
too numerous): they are finite-dimensional. 
On the other hand, the orbits are not trivial, and they increase in
size as the uniton number increases.  The dressing action may therefore
be used to \ll move solutions around\rrr, and it turns out that this is
sufficient to obtain some basic global results on spaces of harmonic maps.

In the general context of zero curvature equations, the dressing
action is defined using a factorization of matrix-valued functions
(the type of factorization depending on the particular problem). It
is usually very difficult to carry out such factorizations
explicitly.  In the case of harmonic maps of finite uniton number, a
drastic simplification is possible (see \cite{Gu-Oh}): the dressing
action on harmonic maps $\Si\to G$ is equivalent to the natural action
of the complex loop group $\LapGC$ on the target space $\Om G=\LaGC/\LapGC$
of the corresponding extended solutions $\F:\Si\to \Om G$.  In terms
of complex extended solutions, this is the action defined 
earlier in Proposition 1.7 (ii).  
If we restrict attention to harmonic maps of uniton number at most $k$, then
the dressing action reduces, roughly speaking, to the action of elements of
$\LapGC$ of the form $\sum_{i=0}^k A_i\la^i$ (assuming a suitable matrix
representation of $G$). This shows that the dressing orbits are finite-dimensional.

Further reductions are possible in the case of special types of maps.
For example, in the case of harmonic maps
associated to $S^1$-invariant extended solutions,
the dressing action reduces to the natural action of the finite-dimensional
Lie group $\GC$ on the \ll twistor space\rr $\Om$. 
In \cite{Gu-Oh}, this concrete realization of the dressing action was
used to enumerate the connected components of the space $\Harm(S^2,G/K)$
of harmonic maps, for $G/K=S^n$ and $G/K=\C P^n$.  

The method involves continuously deforming a general harmonic map,
through a family of harmonic maps, to one of simpler type. In each case a suitable
deformation is obtained by choosing a one parameter subgroup
of $\GC$ and applying it (via the dressing action) to an extended solution. In order to
understand the effects of such deformations, it is very helpful to interpret
the action of a one parameter subgroup of $\GC$ on the generalized
flag manifold $\Om$ as the flow of a Morse-Bott function.

A similar method was used in \cite{Fu-Gu-Ko-Oh} to compute the
fundamental group of $\Harm(S^2,S^n)$. Results on connected components
of $S^1$-invariant extended solutions for other $G/K$ have
been given in \cite{Mu1}, \cite{Mu2}.  (The special feature of the cases $G/K=S^n$ or $\C
P^n$, however, is that {\it all} harmonic maps $S^2\to G/K$ correspond to 
$S^1$-invariant extended solutions.)  Dong (\cite{Dn}) used this method to obtain an
estimate relating the energy and the uniton number of a harmonic map
$S^2\to U_n$.

For harmonic maps which do
not necessarily correspond to  $S^1$-invariant extended solutions,
the same method of Morse-theoretic deformations is available, working
now in the loop group $\Om G$, but these do not seem powerful enough to
determine the connected components of $\Harm(S^2,G/K)$ or $\Harm(S^2,G)$.
We shall give a method in (4) for the case 
$\Harm(S^2,\grkcn)$ or $\Harm(S^2,U_n)$, which makes
use of the canonical form for complex extended solutions.

Finally we note that the space of harmonic maps of finite uniton
number and fixed energy,  which (from the canonical form) is an algebraic variety,
has singularities in general.  However, for $n=1,2$,
(each connected component of) the space 
$\Harm(S^2,\C P^n)$ is actually a smooth manifold. For $n=1$ this is
elementary, and for $n=2$ the manifold structure was established and studied
in \cite{Cr}, \cite{Le-Wo1}, \cite{Le-Wo2}.

\no{\it (4) Topological properties of complex extended solutions}

As in (2) we shall restrict attention to the group $G=U_n$, in order
to make use of the Grassmannian model $\Grn=\Om U_n$ and the
corresponding version $\la W^\pr\sub W$ of the harmonic map equation.
The great advantage of this formulation, for harmonic maps of
finite uniton number, is that $W:\Si\to Gr_r(\Hn_+/\la^k \Hn_+)$ is a
holomorphic map to a finite-dimensional Grassmannian.  We shall
also restrict attention to the case $\Si=S^2$.   The homotopy class $\vert W\vert$ of
$W$ is then a non-negative integer, which is known to be (a constant multiple
of) the energy of the corresponding harmonic map.  We shall use complex
extended solutions to prove Theorem 2.9 below, that any two harmonic maps with the
same energy may be connected by a continuous path in the space of
harmonic maps.

We begin with some remarks on general holomorphic maps $V:S^2\to\grkcn$.
Using the Schubert cell decomposition of $\grkcn$, and the fact that $S^2$ is
one-dimensional, we may write $V=\H E$ for some $k$-plane $E$,  where 
$\H:S^2-F\to \glnc$ is a holomorphic map from the complement of a finite set $F$.
Moreover, we may write $\H=\exp\,B$ where 
$B:S^2-F\to \n$ is a holomorphic map
into some nilpotent Lie subalgebra $\n$ of $\frakglnc$. 
(The affine space $(\exp\,\n)E$ is a Schubert cell
in $\grkcn$, biholomorphically equivalent to $\n$ itself,
and the image of $V\vert_{S^2-F}$
lies in this cell.)   The map $B$ is globally defined on
$S^2$ if and only if it (and hence $V$) is constant.  How, then, does (the locally
defined) $B$  reflect the (global) invariant $\vert V\vert$ ? The answer is provided
by classical Schubert calculus:

\proclaim{Proposition 2.5} Let $Z$ be a linear subspace of $\C^n$ of
codimension $k$. Then $\vert V\vert$ is the number of points $z\in S^2$ 
(counted with multiplicities) such that 
$\dim V(z)\cap Z \ge 1$, whenever this number is finite.
\endproclaim

\no If there exist points $z_1,\dots,z_d$ in the domain $S^2-F$ of $B$, such
that $d=\vert V\vert$ and $\dim (\H(z_i)E)\cap Z \ge 1$ for all $i$,
then we must have $\{z_1,\dots,z_d\}=\{z\in S^2\st \dim V(z)\cap Z \ge 1\}$.
By general position arguments, it is clear that, for
any $B$, there exist $Z$ and $z_1,\dots,z_d$ with these properties.
In this sense, $B$ carries sufficient information to determine $\vert V\vert$.
We shall use later on a partial converse of this statement.  Namely, for
any $V$ and any $Z$, it is possible to deform $V$ continuously to a
holomorphic map $\tilde V=\tilde\H E$ such that there exist 
$\tilde z_1,\dots,\tilde z_d$ with
$\{\tilde z_1,\dots,\tilde z_d\}=\{z\in S^2\st \dim \tilde V(z)\cap Z \ge 1\}$,
and $\{\tilde z_1,\dots,\tilde z_d\} \sub S^2 - F$.

\no{\it Example 2.6:} In the case $k=1$ and 
$E=\{ (x,0,\dots,0) \st x\in\C \} = \C$, a
holomorphic map $V=\H E$ is the same thing as a map
$V:S^2\to\C P^{n-1}$ which can be written in homogeneous coordinates
as $V(z)=[1;a_1(z);\dots; a_{n-1}(z)]$, for some  rational
functions $a_1,\dots,a_{n-1}$.  In this case the finite set $F$
is the union of the poles of $a_1,\dots,a_{n-1}$.
Let us choose the codimension one subspace
$Z=\{ (x_1,\dots,x_{n-1},0) \st x_1,\dots,x_{n-1}\in\C \}$. Then,
assuming that $a_{n-1}$ is not identically zero,
the homotopy class of $V$ is measured by the zeros of $a_{n-1}$ together
with certain poles of $a_1,\dots,a_{n-1}$. By a continuous
deformation of $a_{n-1}$ --- and this applies even in the case
where $a_{n-1}$ is identically zero --- we may deform $V$
to $\tilde V$ in such a way that the homotopy class of $\tilde V$ is 
measured exactly by the zeros of $\tilde a_{n-1}$.\qed

Now we return to the case of a complex extended solution $W$. Let
$$
W=\H \ga_v \Hn_+ = (\exp\,C)\ga_v\Hn_+
= (A_0+\la A_1 + \dots + \la^{k-1} A_{k-1})\ga_v
\Hn_+ \quad\text{(say)}
$$
be the canonical form from (1). We have
$$
\ga_v\Hn_+=E_0 \oplus \la (E_0\oplus E_1) 
\oplus \dots \oplus\la^{k-1}(E_0\oplus\dots\oplus E_{k-1})
\oplus \la^k \Hn_+
$$
where $\C^n=E_k\oplus E_{k-1}\oplus\dots \oplus E_0$ is the common eigenspace
decomposition of the matrices $\ga_v(\la)$ (for all $\la$);
thus $\dim E_i=a_i$. Let
$$
X_0=H E_0 = (A_0+\la A_1 + \dots + \la^{k-1} A_{k-1})E_0.
$$
This defines a holomorphic map $X_0:S^2\to Gr_{a_0}(\C^{kn})$,
which can be taken as the first of a sequence $X_0,X_1,\dots,X_{k-1}$ 
producing $W$ in the manner of Theorem 2.4.

Since $X_0$ is part of the unconstrained holomorphic data
defining $W$, any deformation of $X_0$ through holomorphic maps
gives rise to a deformation of $W$ through extended solutions.
In other words, given a family $X_0^t$ with $0\le t\le 1$,
we define
$W_t=X_t+\la X_t^\pr + \la^2 X_t^{\pr\pr} + \dots + 
\la^{k-1} X_t^{(k-1)} + \la^k\Hn_+$, where
$X_t=X_0^t+\la X_1 + \la^2 X_2 + \dots + \la^{k-1} X_{k-1}$.
(It turns out that we shall not need to deform $X_1,\dots,X_{k-1}$.)
Our objective is to find a deformation $X_0^t$ which gives rise
to a {\it continuous} deformation $W_t$ of $W$, with the additional
property that $W_1$ has lower uniton number than $W_0 (=W)$.
We shall obtain such a deformation by deforming continuously
the (coefficients of the) component functions of $\H$.  To verify
the continuity, we need the following proposition and its proof.

\proclaim{Proposition 2.7} $\vert W\vert = \vert X_0\vert$.
\endproclaim

\no A proof is given in Appendix C.

We can use this in the following way
to determine when the deformation $X_0^t$
gives rise to a continuous deformation of $W$.  Since we consider only
deformations obtained by changing the coefficient functions
continuously, any discontinuity will force $\vert W_1\vert < \vert W_0\vert$.
(This may be seen by regarding  $W_0$ as a holomorphic
map $S^2\to\C P^N$ for some large $N$; for maps into projective
space the assertion is clear.)  Therefore, if $\vert W_t\vert$ remains constant,
the deformation must be continuous.   Let $z_1,\dots,z_d$ be
the points where the first entry of the last column of $A_{k-1}$
vanishes (without loss of generality we may assume these points are in the domain
of $\H$, and, by the proof of Proposition 2.7, 
we may assume $d=\vert W\vert = \vert X_0\vert$).  If
these points are unchanged during the deformation, then we must have
$\vert W_t\vert \ge \vert W_0\vert$ for all $t$, hence
$\vert W_t\vert = \vert W\vert$, and hence any deformation of this type is
continuous.  Our strategy will be to deform $X_0$ 
by deforming only $A_0E_0$ (and we will
assume $k\ge 2$, so that $A_0\ne A_{k-1}$).  Then $z_1,\dots,z_d$ will
indeed remain unchanged.  

\proclaim{Proposition 2.8} Assume $k\ge 2$.  Then there exist (a sequence of) continuous
deformations of the above type such that the resulting extended solution has
uniton number less than $k$.
\endproclaim

\no A proof is given in Appendix C.

It follows that all non-constant extended solutions can be deformed
continuously to extended solutions of uniton number $1$.  The
latter correspond to holomorphic maps into Grassmannian submanifolds
of $U_n$, and (using the method of (3) above) it is easy to show that 
any two such may be connected by a continuous path, if they have
the same energy.
Let $\Harm_e(S^2,U_n)$ be the space of harmonic maps from $S^2$ to $U_n$
which have energy $e$. (We may assume that the energy has been normalized so
that $e=\vert W\vert$, where $W$ is any corresponding extended solution.)
Then the above method gives:

\proclaim{Theorem 2.9}  For any $e\ge 0$, and any $n\ge 1$,
$\Harm_e(S^2,U_n)$ is path-connected.
\endproclaim

In Appendix C, in addition to the omitted proofs, we shall give a 
much more explicit method for the case $n=3$.
We remark that Anand (\cite{An}) has obtained independently this and other results on
the space $\Harm_e(S^2,U_3)$ for low values of $e$, by using the
vector bundle approach.

\no{\it (5) Harmonic maps of finite uniton number into symmetric spaces}

If $G/K$ is a symmetric space of $G$, harmonic maps $\Si\to G/K$ 
may be considered as special examples of harmonic maps $\Si\to G$,
by the remarks in Appendix A.  To what extent can the results of
(1)-(4) be adapted to such maps?  A partial answer to this question is
that {\it all\,}  the above results can be extended to the case of {\it inner}
symmetric spaces, i.e. where the involution of the symmetric space is
an inner automorphism of $G$.

A canonical form for the case of maps into an inner symmetric space
$G/K$ was given in \cite{Bu-Gu}; this is a straightforward modification
of (1).  By making use of this, the Frenet frame construction of (2)
may be extended to the case of maps into any inner symmetric space
associated with $G=U_n$, i.e. any Grassmannian $G/K=\grkcn$.
The Morse-theoretic deformations of (3) were used primarily in the
case of maps into symmetric spaces in the first place, so no further
comment is needed here. Finally, in the case of (4), one can prove
the analogue of Theorem 2.9 for the space $\Harm(S^2,\grkcn)$, again
by studying the associated complex extended solutions. That is, if
$\Harm_{d,e}(S^2,\grkcn)$ is the space of harmonic maps from $S^2$ to $\grkcn$
which have energy $e$ and homotopy class $d$, then this space is path-connected.
The method of proof of Theorem 2.9 permits a reduction to the case $k=1$,
where the result is given by combining \cite{Cr} and \cite{Gu-Oh}. 

The situation for general (compact) non-inner symmetric spaces 
has not been investigated.
Since any such space admits a totally geodesic Cartan embedding into a
corresponding compact Lie group, the situation is in principle
covered by the results for Lie groups. However, a more
precise analysis would be desirable.

\head
Appendix A: Review of harmonic maps, from Frenet frames to the zero curvature equation
\endhead

These brief historical comments provide the differential geometric
context for the zero curvature equation method, and
go as far as the 1989 articles of Uhlenbeck (\cite{Uh}) and Segal (\cite{Se})
which we take as our starting point in \S 1.  There are at least
two deficiencies in our exposition.  One is that we ignore a vast amount of
exploratory work in differential geometry which, while not central to our story,
was instrumental in shaping attitudes to research on harmonic
maps from Riemann surfaces to Lie groups or homogeneous spaces.   For this
we refer to the comprehensive survey articles \cite{Ee-Le1}, \cite{Ee-Le2}.
The other is that we omit references to previous work in integrable systems  
which provided the foundation for \cite{Uh} and \cite{Se}.  

Harmonic maps from the two-sphere $S^2$
to the complex Grassmannian manifold $\grkcn$ or the unitary group $U_n$ are
\ll model examples\rr for two important problems: (1)  Minimal
immersions of (compact) Riemann surfaces into symmetric spaces
(differential geometry), and (2) The nonlinear $\si$-model or chiral model
(mathematical physics). The case
$G/K=\grkcn$ is closely related to the case $G=U_n$, so we shall concentrate on the latter.
Remarks on general Lie groups $G$ and general symmetric spaces
$G/K$, as well as general Riemann surfaces $\Sigma$ in the
role of the domain of $\phi$,  will appear later on.

The harmonic map equation, for maps $\phi:S^2\to U_n$, is the second order
partial diferential equation
$$
(\phi^{-1} \phi_{\bar z})_z + (\phi^{-1} \phi_{z})_{\bar z} = 0
$$
where $z$ is the usual local coordinate on $S^2=\C\cup\infty=\C P^1$.
It is the Euler-Lagrange equation for the functional $\int \vert\vert
d\phi \vert\vert^2$, and is therefore a very natural generalization of the
geodesic equation in Riemannian geometry.  Harmonic maps $S^2\to\grkcn$
constitute a special case, in the following sense.  The manifold
$\grkcn$ can be embedded totally geodesically in $U_n$ by sending a complex $k$-plane $V$
to the unitary linear transformation $\pi_V - \pi_{V^\pe}$, where $\pi_V:\C^n\to\C^n$
is orthogonal projection onto $V$ (with respect to usual the Hermitian inner product
of $\C^n$), and a map $\phi:S^2\to \grkcn$ is harmonic if and only if
the composition $S^2\to \grkcn\to U_n$ is harmonic.

The above equation admits many solutions, the simplest (non-constant) ones
being the holomorphic maps $S^2\to\grkcn$.   The space $\Hol(S^2,\grkcn)$ of all
such holomorphic maps has connected components $\Hol_d(S^2,\grkcn)$ indexed by
$d=0,1,2,\dots$, where $d$ is the homotopy class $\vert\phi\vert\in\pi_2\grkcn$.
Each connected component has the structure of a  complex
manifold of dimension $nd + k(n-k)$ 
(noncompact if $d\ge 1$), which has been thoroughly
studied from the point of view of algebraic topology.  The case $k=1$, $n=2$ (hence
$\grkcn=S^2$) is paricularly simple, since holomorphic maps $S^2\to S^2$
are just the same as rational functions of the complex variable $z$.

In contrast, the larger class of (not necessarily holomorphic) harmonic
maps is much less easy to describe.  The first major result was for $k=1$, i.e. harmonic
maps $S^2\to \C P^{n-1}$, and in the related case of harmonic maps $S^2\to S^m$.  From the
point of view of differential geometry, this originated from work of
Calabi and of Chern on the Cartan moving frame of a minimal immersion, and was achieved
by several people simultaneously in the late 1970's; the following version is that of
Eells and Wood (\cite{Ee-Wo}). Let $[f]:S^2\to \C P^{n-1}$ be any holomorphic
map, and let $i\in\{0,1,\dots,n-1\}$. Define $\phi:S^2\to \C P^{n-1}$ by
$$
\phi(z)=([f(z)]+[f^\pr(z)]+\dots+[f^{(i)}(z)]) \ominus
([f(z)]+[f^\pr(z)]+\dots+[f^{(i-1)}(z)]).\tag $\ast$
$$
Here we regard $f,f^{\pr},\dots$ as $\C^n$-valued rational functions,
and $[f(z)]+[f^\pr(z)]+\dots$ means the vector space spanned by the lines
$[f(z)],[f^\pr(z)],\dots$. The
notation $A\ominus B$ denotes the (Hermitian) orthogonal complement of
$B$ in $A$.  Formula $(\ast)$ defines a $\C P^{n-1}$-valued function which
is smooth except at a finite number of singular points; it can be shown
that the singularities are removable.  (If any of the $\C^n$-valued functions appearing
here are identically zero, an obvious modification of the construction produces
maps into lower dimensional complex projective spaces.)  The result is:

\proclaim{Theorem A.1} The map $\phi$ defined by formula $(\ast)$ 
is harmonic, and all harmonic maps  $S^2\to \C P^{n-1}$ arise this way.
\endproclaim

This is a very satisfactory solution to the problem: all harmonic maps
$S^2\to \C P^{n-1}$ are obtained by performing a series of algebraic
operations and derivatives on rational functions.  Rational functions
are \ll known\rr objects and the operations are mechanical (in particular,
programmable on a computer), so the solution is as explicit as one
could hope for. The computation of the orthogonal complement is somewhat
messy, of course, and so the formula for the homogeneous components
of $\phi$ (in terms of the components of $f$) may be
unpleasant to look at.  But this is not a defect of the method; on the contrary,
it tells us that the formula $(\ast)$ is the heart of the matter.

The simplicity of this result has been both a blessing and a curse.
It inspired great efforts to find generalizations to other situations,
and it demonstrated important themes in mathematical physics such as twistor theory
and the role of instanton solutions of field equations.
On the other hand, generalizations to the case of other target spaces
were frustratingly slow in coming. For several important
manifolds $M$, such as Grassmannians, quadrics, real and quaternionic projective
spaces,  various results pointed to a description of all harmonic maps
$S^2\to M$ in terms of \ll holomorphic data\rrr, but 
always lacked the simplicity of the $\C P^{n-1}$ case.
Surveys of these purely differential geometric developments
can be found in \cite{Ee-Le1}, \cite{Ee-Le2}.

The next major progress came from the theory of integrable systems.
It had been known for some time that the harmonic map equation, like several other
important equations, could be written as a zero curvature equation.
To do this, one introduces the matrix-valued complex $1$-form 
$$
\om = \frac12 (1-\frac1{\la}) \phi^{-1} \phi_{z} dz
+\frac12 (1-{\la}) \phi^{-1} \phi_{\bar z} d\bar z
$$
where $\la$ is a complex parameter.  Then an easy computation
shows that the harmonic map equation is equivalent to the condition
$$
d\om + \om \wedge \om = 0\quad\text{for all}\quad \la.
$$
(This is the condition that the curvature
tensor of the connection $d+\om$ is identically zero.)
We refer to \cite{Uh}, \cite{Se}, \cite{Hi-Se-Wa} for background
information on this somewhat mysterious formulation.

From this point of view, the first significant advance was made
by Uhlenbeck in \cite{Uh}, where the concept of \ll extended solution\rr
was introduced.  An extended solution associated to a harmonic map
$\phi:S^2\to U_n$ is any map $\F:S^2\times\C^\ast\to U_n$ such that $\om = \F^{-1}d\F$.  
The existence of such a map $\F$,
smooth in the variable $z\in S^2$ and holomorphic in the variable
$\la\in \C^\ast = \C-\{0\}$,  is an elementary
fact (see, for example,  \cite{Sh}). Conversely, if $\F$
is a map with the properties
$$
\align
\F^{-1} \F_{z} &= \text{a function of $z$ and $\la$ which is linear in $\frac 1\la$} \\
\F^{-1} \F_{\bar z} &= \text{a function of $z$ and $\la$ which is linear in $\la$ }
\endalign
$$
then the map $\phi(z)=\F(z,-1)$ satisfies the zero curvature equation 
and is therefore harmonic.

This reformulation of the harmonic map equation has two specific advantages. 
First, by restricting $\F$ to the circle $\vert \la\vert = 1$, we may regard $\F$
as a map from $S^2$ into the (based) loop group $\Om U_n = \{ \ga:S^1\to U_n
\st \ga(1)=I \}$, and the second condition implies that this
$\F:S^2\to \Om U_n$ is holomorphic with respect to the standard complex
structures (see \cite{Pr-Se} for the complex structure of $\Om U_n$).
Thus, the harmonic condition is immediately translated to holomorphic one,
and from now on one works entirely with holomorphic objects.
The second advantage is based on the (at first sight inconvenient) fact that the
correspondence between $\phi$ and $\F$ is not one to one.  If $\F$ is an extended
solution associated to $\phi$, then $\ga \F$ is obviously an extended solution
associated to $\ga(-1) \phi$, and it is another
elementary fact that this is essentially the only ambiguity.  
The freedom to pre-multiply $\F$ by a loop $\ga$ may be used to
seek a canonical representative associated to a given harmonic map $\phi$.
Uhlenbeck used this freedom to prove (in \cite{Uh}) the following statement:

\proclaim{Theorem A.2} For any harmonic map
$\phi:S^2\to U_n$, there exists an associated extended solution $\F$ of the form
$\F(z,\la)= \sum_{i=0}^k A_i(z) \la^i$, with $k\le n-1$.
\endproclaim

\no The least such (non-negative) integer $k$ is
called the (minimal) uniton number of $\phi$. This is a measure
of the complexity of $\phi$; for example, constant maps have uniton number
zero, and holomorphic maps $S^2\to\grkcn(\to U_n)$ have uniton
number one.

A second result of \cite{Uh} was the existence of a multiplicative \ll uniton
factorization\rr of $\F$ into 
$k$ factors of the form $\pi_{V(z)} + \la \pi_{V(z)^\pe}$, where
each $V$ is a (smooth) map into a complex Grassmannian.  
When $k=1$ there is only one factor
and the map $V$ is holomorphic. When $k\ge 2$ the maps $V$ are in some sense
holomorphic with respect to perturbations of the standard complex
structure, where the perturbation depends on
the previous factors.  Roughly speaking, this means that a
harmonic map of uniton number $k$ can be factored as a product of
$k$ holomorphic maps (but, because of
the perturbations, this statement is not strictly true, which
reduces the usefulness of the factorization as a practical tool).

A heuristic explanation of this factorization was given in \cite{Uh}, by using the
\ll dressing\rr procedure of integrable systems theory.  This procedure
generally works for zero curvature equations which involve a condition on the
type and location of the poles of $\om$.  
(In the case of the harmonic map equation, there
are two simple poles, at $\la=0$ and $\la=\infty$.) Given a solution $\F$,
one obtains under certain conditions a new solution of the form
$\F X$, where $X$ is the solution of an associated \ll Riemann-Hilbert\rr
problem for a contour consisting of two small circles around the poles.  We omit
further details as we use a different manifestation of the dressing procedure 
in this article; the version of \cite{Uh} and its relation to the uniton factorization
is discussed in more detail in \cite{Be-Gu}.

The treatment of Uhlenbeck is analytic, deriving finiteness of the
uniton number from ellipticity of the harmonic map equation and
compactness of $S^2$.   
A different version was given by Segal in \cite{Se}, using the
Grassmannian model of the loop group $\Om U_n$; here the finiteness of the uniton number
is a consequence of having a holomorphic map from the compact manifold $S^2$
into $\Om U_n$. Without going into details, we remark that the loop group
$\Om U_n$ may be identified with an infinite-dimensional Grassmannian
manifold $\Grn$, which is a subspace of the Grassmannian of all
linear subspaces of the Hilbert space
$\Hn=L^2(S^1,\C^n)$ (see \cite{Pr-Se}),  and the conditions that a
(smooth) map
$W:S^2\to \Grn$ corresponds to an extended solution $\F:S^2\to\Om U_n$ are
$$
\align
W_z&\sub\frac1\la W\\
W_{\bar z}&\sub W.
\endalign
$$
The first condition means that the vector $\b s(z)/\b z$ is
contained in the subspace $\la^{-1} W(z)$ of $\Hn$, for every 
(smooth) map $s:S^2\to \Hn$ such that $s(z)\in W(z)$.
The second condition is interpreted in a similar way; it is equivalent
to saying that the map $W$ is holomorphic.  Segal's version of the
finiteness of the uniton number is that the map $W$ (corresponding
to the original harmonic map $\phi:S^2\to U_n$) may be chosen so that
its image lies in the \ll algebraic\rr Grassmannian
$\Grnalg$.  An elementary treatment of the loop group
and Grassmannian model approach to harmonic maps may be found in \cite{Gu}.

Finally, a brief comment on the situation for maps $\Si\to G$ (or $G/K$),
where $\Si$ is an arbitrary compact Riemann surface and $G$ is an
arbitrary compact Lie group (or $G/K$ a symmetric space).  The harmonic map
equation makes sense here, with $z$ interpreted as a local coordinate
on $\Si$, and $f^{-1}df$ interpreted as the pullback
by $f$ of the Maurer-Cartan form on $G$, as does the definition of
extended solution.  An extended solution gives rise to a harmonic
map in the above manner, {\it but there is no
guarantee that a harmonic map $\Si\to G$ gives rise to an extended solution} 
$\Si\to \Om G$, when $\Si$ is not simply connected 
(i.e. when $\Si$ has positive genus).  The class of harmonic maps 
which admit extended solutions $\Si\to \Om G$ is precisely the class
of harmonic maps of finite uniton number.

\head
Appendix B: Proof of Theorem 2.4
\endhead

We begin with extended solutions of type $(2,1,0)$.  From
Example 1.10 we have the canonical form 
$W(z,\la=\exp\, C(z,\la) \diag(\la^2,\la,1) \Hn_+$
where $n=3$ and
$$
C
=
\pmatrix
0 & a & b\\
0 & 0 & c\\
0 & 0 & 0
\endpmatrix
+
\la
\pmatrix
0 & 0 & d\\
0 & 0 & 0 \\
0 & 0 & 0 
\endpmatrix.
$$
After making the change of variables in Example 2.1 we obtain
$$
W = \exp\,C\diag(\la^2,\la,1)\Hn_+
=
\left[
\pmatrix
1 & \ga & \al\\
0 & 1 & \be\\
0 & 0 & 1
\endpmatrix
+
\la
\pmatrix
0 & 0 & \de\\
0 & 0 & 0 \\
0 & 0 & 0 
\endpmatrix
\right]
\diag(\la^2,\la,1)
\Hn_+.
$$
Now, $\diag(\la^2,\la,1) \Hn_+=V_3 \oplus \la V_3 \oplus \la V_2 \oplus \la^2\Hn_+$,
where $\C^3=V_1\oplus V_2\oplus V_3$ is the decomposition given by the directions of the
standard basis vectors. If we introduce
$$
l=
\pmatrix
\al\\
\be\\
1
\endpmatrix,
m=
\pmatrix
\de\\
{0}\\
{0}
\endpmatrix,
n=
\pmatrix
\ga\\
1\\
{0}
\endpmatrix,
$$
then $W$ is the span of the three vector-valued functions
$l+\la m$, $\la(l+\la m)$, $\la n$ (mod $\la^2\Hn_+$). 
The equation $\ga=\al^\pr/\be^\pr$ (the complex extended solution equation)
implies that $[l^\pr]=[n]$. Hence:
$$
\align
W&=[l+\la m]+\la[l+\la m]+\la[l^\pr]+\la^2\Hn_+\\
&=[l+\la m]+\la[l+\la m]+\la[l^\pr+\la m^\pr]+\la^2\Hn_+\\
&=[l+\la m]+\la [l+\la m]^\pr+\la^2\Hn_+.
\endalign
$$  
This establishes the formula of Theorem 2.4 for harmonic maps of type $(2,1,0)$.
For types $(1,1,0)$ and $(1,0,0)$ the statements of Theorem 2.4 are obvious.

In the case of $U_4$, the arguments for type $(3,2,1,0)$ and $(2,2,1,0)$
are exactly like that for $(2,1,0)$ above, so let us consider type $(2,1,1,0)$.
We may write
$$
\exp\,C=
\pmatrix
1 & \ga & \be & \al_1 \\
0 & 1 & 0 & \al_2\\
0 & 0 & 1 & \al_3\\
0 & 0 & 0 & 1
\endpmatrix
+\la
\pmatrix
0 & 0 & 0 & \de\\
0 & 0 & 0 & 0\\
0 & 0 & 0 & 0\\
0 & 0 & 0 & 0
\endpmatrix
$$
(thus defining the change of variable for this situation),
in terms of which the complex extended solution equation is
$$
\al_1^\pr=\ga\al_2^\pr + \be\al_3^\pr.
$$
Let us define
$$
l=
\pmatrix \al_1\\ \al_2\\ \al_3\\ 1 \endpmatrix,\quad
m=
\pmatrix \de \\ 0 \\ 0 \\ 0 \endpmatrix,\quad
p=
\pmatrix \be \\ 0 \\ 1 \\ 0 \endpmatrix,\quad
q=
\pmatrix \ga \\ 1 \\ 0 \\ 0 \endpmatrix.
$$
Since 
$\diag(\la^2,\la,\la,1)\Hn_+=V_4\oplus\la V_4 
\oplus \la V_3 \oplus \la V_2 \oplus \la^2\Hn_+$,
we have
$$
W=[l+\la m] + \la[l+\la m] + \la [p] + \la [q] + \la^2\Hn_+.
$$
The complex extended solution equation implies
$[l^\pr]\sub [p]+[q]$, so
$$
W=[l+\la m] + \la[l+\la m]^\pr + \la[n] + \la^2\Hn_+
$$
where $[n]$ is any complementary (one-dimensional) subspace to
$[l^\pr]$ in $[p]+[q]$. If we put $X_0= [l+\la m]$ and $X_1=[n]$,
and $X=X_0+\la X_1$, then $W=X+\la X^\pr + \la^2\Hn_+$,
as required.

The next case to consider is $(2,1,0,0)$. We begin as usual with
the canonical form, which is
$$
\exp\,C=
\pmatrix
1 & \ga & \be_1 & \al_1 \\
0 & 1 & \be_2 & \al_2\\
0 & 0 & 1 & 0\\
0 & 0 & 0 & 1
\endpmatrix
+\la
\pmatrix
0 & 0 & \ep & \de\\
0 & 0 & 0 & 0\\
0 & 0 & 0 & 0\\
0 & 0 & 0 & 0
\endpmatrix.
$$
Then we express $W$ in terms of 
$$
l_1=
\pmatrix \al_1\\ \al_2\\ 0\\ 1\endpmatrix,\quad
l_2=
\pmatrix \de\\ 0\\ 0\\ 0\endpmatrix,\quad
m_1=
\pmatrix \be_1\\ \be_2\\ 1\\ 0\endpmatrix,\quad
m_2=
\pmatrix \ep\\ 0\\ 0\\ 0\endpmatrix,\quad
r=
\pmatrix \ga\\ 1\\ 0\\ 0\endpmatrix,
$$
that is
$$
W=[l_1+\la l_2] + \la [l_1+\la l_2] + [m_1+\la m_2] + \la [m_1+\la m_2] 
+ \la[r] + \la^2\Hn_+.
$$
The complex extended solution equations are
$\al_1^\pr = \de \al_2^\pr$, $\be_1^\pr = \de \be_2^\pr$, which
imply $[l_1^\pr]=[m_1^\pr]=[r]$. 

At this point we do not appear to have unconstrained holomorphic
data, because the four vector-valued functions $l_1,l_2,m_1,m_2$
are required to satisfy the condition $[l_1^\pr]=[m_1^\pr]$.  There
are two ways of dealing with this.  The first is to observe that
$X=X_0=[l_1+\la l_2] + [m_1+\la m_2]$ always generates a solution
$W=X+\la X^\pr+\la^2\Hn_+$ of the equation $\la W^\pr\sub W$,
which coincides with one of type $(2,1,0,0)$ in the particular case
$[l_1^\pr]=[m_1^\pr]$. (In the other cases, a dressing transformation converts
it to a solution of simpler type.) The second point of view is
to use the fact that two vector-valued functions $l_1,m_1$ satisfy
$\dim [l_1]^\pr + [m_1]^\pr = 3$ if and only if there exists some
$l$ such that $[l]\sub[l_1]+[m_1]$, $[l]^\pr=[l_1]+[m_1]$,  
and $[l]^{\pr\pr}=[l_1]^\pr +[m_1]^\pr$. 
Using this we can rewrite the above formula for $W$ as
$$
\align
W
&=
[l_1+\la l_2] + \la [l_1+\la l_2] + [m_1+\la m_2] + \la [m_1+\la m_2] 
+ \la[r] + \la^2\Hn_+\\
&=
[l_1+\la l_2] + \la [l_1] + [m_1+\la m_2] + \la [m_1] 
+ \la[r] + \la^2\Hn_+\\
&=
[l_1+\la l_2] + [m_1+\la m_2] + \la ([l_1]+[m_1] +[r]) + \la^2\Hn_+\\
&=
[al+bl^\pr+\la l_2] + [cl+dl^\pr+\la m_2] + \la ([l]+[l^\pr] +[l^{\pr\pr}])
 + \la^2\Hn_+
\endalign
$$
(writing $l_1=al+bl^\pr$ and $m_1=cl+dl^\pr$ where $a,b,c,d$ are
meromorphic functions).  This can be written in the form
$[l+\la m] + [l^\pr+\la n] + \la ([l]+[l^\pr] +[l^{\pr\pr}]) + \la^2\Hn_+$
as stated in Theorem 2.4.

For the remaining types $(1,1,1,0)$, $(1,1,0,0)$ and $(1,0,0,0)$,
the statements of Theorem 2.4 are obvious.

\head
Appendix C: Proof of Propositions 2.7 and 2.8
\endhead

\demo{Proof of Proposition 2.7} 
We regard $W$ as a holomorphic map $S^2\to Gr_r(\C^{kn})$,
where
$$
\align
r&= \dim \ E_0 \oplus \la(E_0\oplus E_1) \oplus \dots 
\oplus\la^{k-1}(E_0\oplus\dots\oplus E_{k-1})\\
&= a_0 + (a_0+a_1) + \dots + (a_0+\dots+a_{k-1}).
\endalign
$$
To measure $\vert W\vert$ (using Proposition 2.5) we need a linear
subspace $Z$ of codimension $r$ in $\C^{kn}$.  We choose
$$
Z=
(E_k\oplus\dots\oplus E_1)\oplus V_n\oplus
\la(E_k\oplus\dots\oplus E_2)\oplus
\dots\oplus
\la^{k-2}(E_k\oplus E_{k-1})\oplus
\la^{k-1}(E_k\ominus V_1)
$$
where (as usual) $V_1,\dots,V_n$ are the standard basis directions in $\C^n$
(so $\C^n=V_1\oplus\dots \oplus V_n$ and $V_n\sub E_0$, $V_1\sub E_k$).
On the other hand, $X_0$ defines a holomorphic map $S^2\to Gr_{a_0}(\C^{kn})$,
and in fact a holomorphic map into the Grassmannian of $a_0$-dimensional
subspaces of 
$$
\C^n\oplus
\la (E_k\oplus\dots\oplus E_2)\oplus
\la^2 (E_k\oplus\dots\oplus E_3)\oplus
\dots\oplus
\la^{k-1}E_k.
$$
Now, $Z$ is a subspace of this as well, and it has codimension $a_0$.
So $Z$ may be used to measure $\vert X_0\vert$.
The proof of the proposition is completed by noting that, by construction,
$\dim W(z)\cap Z\ge 1$ if and only if $\dim X_0(z)\cap Z\ge 1$,
since both conditions are equivalent to the vanishing of the first entry
of the last column of  the matrix $A_{k-1}$.
\qed\enddemo

\demo{Proof of Proposition 2.8} We take 
a deformation of the form described in (4) of \S 2 which renders
$A_0E_0$ constant.  For the resulting extended solution we claim that

\no(1)  $\la^k \Hn_+ \oplus \la^{k-1}A_0E_0 \sub W$, and

\no(2) $W \sub \Hn_+ \ominus (A_0E_0)^\perp$

\no(where the orthogonal complement is taken in $\C^n$).
To establish these, we shall use the formula
$$
W=(A_0+\dots +\la^{k-1} A_{k-1})
(E_0 \oplus \la (E_0\oplus E_1) \oplus \dots 
\oplus\la^{k-1}(E_0\oplus\dots\oplus E_{k-1})
\oplus \la^k \Hn_+).
$$
It follows that $W$ contains 
$$
(A_0+\dots +\la^{k-1} A_{k-1})\,
(\la^{k-1}(E_0\oplus\dots\oplus E_{k-1})\oplus \la^k \Hn_+)
=
\la^{k-1}A_0(E_0\oplus\dots\oplus E_{k-1})\oplus \la^k \Hn_+.
$$
Since $A_0E_0$ is constant, we obtain (1). Similarly,
by considering $W+\la \Hn_+$, we see that $W$ is orthogonal
to $(A_0E_0)^\perp$, and this is (2).
Any $W$ satisfying (1) and (2) must have uniton number
less than $k$ (after translation by a constant loop), since
$\dim [\Hn_+\ominus (A_0 E_0)^\perp] / 
[\la^k \Hn_+ \oplus \la^{k-1} A_0 E_0] = (k-1)n$.
\qed\enddemo  

As a concrete illustration we return to
the case $n=3$, where the underlying geometry is more transparent.
We have seen in (1) of \S 2 that there are three types of non-constant
harmonic map $S^2\to U_3$, namely $(2,1,0)$, $(1,1,0)$, and $(1,0,0)$.
The last two correspond to holomorphic maps from $S^2$ into totally
geodesic submanifolds $Gr_2(\C^3)$ and $\C P^2$ of $U_3$, and
(as stated earlier) all such maps are in the same path-component
of $\Harm_e(S^2,U_3)$.  The proof of Theorem 2.9 is carried out by
showing that any harmonic map of type $(2,1,0)$ (hence with uniton number
at most $2$) may be deformed continuously through harmonic maps
to a harmonic map of uniton number $1$.
Let us see how this works, explicitly.

Without loss of generality we begin with a complex extended solution
of the form
$$
\align
W&= (A_0+\la A_1) \diag(\la^2,\la,1) \Ht_+\\
&=
\left[
\pmatrix
\ \ 1\ \  & \al^\pr/\be^\pr & \ \ \al\ \ \\
0 & 1 & \be\\
0 & 0 & 1
\endpmatrix
+
\la
\pmatrix
0 & 0 & \de\\
0 & 0 & 0 \\
0 & 0 & 0 
\endpmatrix
\right]
V_3 \oplus \la (V_3 \oplus V_2) \oplus \la^2\Ht_+
\endalign
$$
(cf. Appendix B).  Here, $\al,\be,\de$ are arbitrary rational functions,
and the zeros of $\de$ are certain points $z_1,\dots,z_e$, where
$e=\vert W\vert$.  These are the points $z$ such that
$\dim W(z)\cap Z \ge 1$, with $Z=\C^3$.  (If $\de$ has less than $e$
zeros in the domain of $A_0+\la A_1$, or if $\de$ is identically
zero, then we deform $\de$ continuously until
we are in the above situation.) 
We have $E_2=V_1$, $E_1=V_2$, $E_0=V_3$, and so
$$
X_0=(A_0+\la A_1)E_0=  
\left[
\pmatrix \al\\ \be \\ 1
\endpmatrix
+
\la
\pmatrix \de\\ 0 \\ 0
\endpmatrix
\right],
\quad\quad
A_0E_0=
\left[ 
\pmatrix \al\\ \be \\ 1
\endpmatrix
\right].
$$

Next, we deform the coefficients of $\al$, $\be$ to make $A_0E_0$
constant. Evidently this is a deformation of $W$ through extended
solutions, and it is a continuous deformation as the zeros of $\de$ are
unchanged.  The resulting complex extended solution, which we may as
well arrange to be
$$
W=
\left[
\pmatrix
1 & 0 & 0\\
0 & 1 & 0\\
0 & 0 & 1
\endpmatrix
+
\la
\pmatrix
0 & 0 & \de\\
0 & 0 & 0 \\
0 & 0 & 0 
\endpmatrix
\right]
V_3 \oplus \la (V_3 \oplus V_2) \oplus \la^2\Ht_+,
$$
satisfies $\la(V_2\oplus V_3) \oplus \la^2\Ht_+ \sub W \sub V_3\oplus\la \Ht_+$, i.e.
$$
\diag(\la^2,\la,\la) \Ht_+ \sub W \sub \diag(\la,\la,1) \Ht_+.
$$
Pre-multiplying by $\diag(\la^{-1},1,1) \Ht_+$ (a trivial dressing
transformation, which does not change the energy of the harmonic
map) shows that
$$
\la \Ht_+ \sub \diag(\la^{-1},1,1) W \sub \diag(1,\la,1) \Ht_+ \sub \Ht_+,
$$
hence $\diag(\la^{-1},1,1) W$ has uniton number $1$.  In fact, $W$
is a very special kind of \ll$1$-uniton\rr as it corresponds to
a holomorphic map $S^2\to S^2$.

So far we have merely repeated the above proof of Theorem 2.9 in
this particular case. A closer examination of the Grassmanian
model (the details of which we leave to the interested reader)
reveals a more geometrical interpretation of this deformation.
The original extended solution 
$W:S^2\to Gr_3( \Ht_+/\la^2\Ht_+) = Gr_3(\C^6)$
actually maps into 
$$
\align
X&=\{ \text{all linear subspaces}\ E\sub\Ht \ \text{such that} \\
&\quad\quad 
\la^2\Ht_+\sub E\sub \Ht_+, \la E\sub E, \dim \Ht_+/E = 3,
\dim (\la^{-1} E\cap \Ht_+)/E\ge 2 \}.
\endalign
$$
This is a four-dimensional algebraic subvariety of $Gr_3(\C^6)$
with one singular point, $E=\la \Ht_+$.  Moreover, $X-\{\la \Ht_+\}$
has the structure of a holomorphic line bundle over the flag manifold
$F_{1,2}(\C^3) = U_3/( U_1\times U_1\times U_1)$ (it is the unstable
bundle of the conjugacy class of the geodesic $\diag(\la^2,\la,1)$
in the Morse theory decomposition of $\Om^{\alg} U_3$, and its
closure --- i.e. $X$ --- is the corresponding Schubert variety in
this algebraic loop group).  To be precise, it is the bundle
$\Cal E_1^\ast \otimes (\C^3/\Cal E_2)$, where $\Cal E_1, \Cal E_2$
are the tautological vector bundles of ranks $1,2$ on 
$F_{1,2}(\C^3)$.  

It follows from this description of the structure of $X$ that the
homotopy class $\vert W\vert$ decomposes as $\vert W\vert = d + d_1 + d_2$,
where

\no(1) $d$ is the number of points (counting multiplicities) in
$W^{-1}(\la \Ht_+)$, and

\no(2) $d_i=-c_1 W_{\infty}^\ast \Cal E_i$, where $W_{\infty} =
\pi\circ W$ and $\pi:X-\{\la\Ht_+\} \to F_{1,2}(\C^3)$ is the
bundle projection.

\no Explicitly, the map $\pi\circ W$ is represented by
$$
\pmatrix
\ \ 1\ \  & \al^\pr/\be^\pr & \ \ \al\ \ \\
0 & 1 & \be\\
0 & 0 & 1
\endpmatrix
V_3 \oplus \la (V_3 \oplus V_2) \oplus \la^2\Ht_+.
$$
The rational function $\de$ represents a meromorphic section of
$\pi$, with $\vert W\vert$ zeros and $d$ poles.  By deforming the
rational functions $\al$, $\be$ to constants, we deform $W$ into
a (compactified) fibre of $\pi$, and hence to a holomorphic
map $S^2\to S^2$.

\eightpoint \it

\no  Department of Mathematics
\newline
Graduate School of Science
\newline
Tokyo Metropolitan University
\newline
Minami-Ohsawa 1-1, Hachioji-shi
\newline
Tokyo 192-0397, Japan

\no martin\@math.metro-u.ac.jp

\enddocument